\documentclass[a4paper,11pt]{amsart}
\usepackage{graphicx}
 \usepackage{mathptmx}
\usepackage{amsmath}
\usepackage{amssymb}
\usepackage{enumitem}
\usepackage{xcolor}

\newmuskip\pFqmuskip

\newcommand*\pFq[6][8]{%
  \begingroup 
  \pFqmuskip=#1mu\relax
  \mathcode`=\string"8000
  \begingroup\lccode`\~=`\,
  \lowercase{\endgroup\let~}\pFqcomma
  F^{#2}_{#3}{\left(\genfrac..{0pt}{}{#4}{#5}\bigg|#6\right)}%
  \endgroup
}
\newcommand{\pFqcomma}{\mskip\pFqmuskip}

\newtheorem{theorem}{Theorem}[section]

\begin{document}

\title[Probabilistic Lah numbers and Lah-Bell polynomials]{Probabilistic Lah numbers and Lah-Bell polynomials }

\author{Yuankui Ma*}
\address{School of Science, Xi’an Technological University, Xi’an 710021, Shaanxi, China}
\email{mayuankui@xatu.edu.cn}

\author{Taekyun  Kim*}
\address{School of Science, Xi’an Technological University, Xi’an 710021, Shaanxi, China;
Department of Mathematics, Kwangwoon University, Seoul 139-701, Republic of Korea}
\email{kimtk2015@gmail.com}
\author{Dae San  Kim* }
\address{Department of Mathematics, Sogang University, Seoul 121-742, Republic of Korea}
\email{dskim@sogang.ac.kr}

\thanks{ * corresponding authors}
\begin{abstract}
Let Y be a random variable whose moment generating function exists in some neighborhood of the origin. The aim of this paper is to study the probabilistic Lah numbers associated with $Y$ and the probabilistic Lah-Bell polynomials associated with $Y$, as probabilistic versions of the Lah numbers and the Lah-Bell polynomials, respectively. We derive some properties, explicit expressions, recurrence relations and certain identities for those numbers and polynomials. In addition, we treat the special cases that $Y$ is the Poisson random variable with parameter $\alpha > 0$ and the Bernoulli random variable with probability of success $p$.
\end{abstract}

\maketitle

\section{Introduction}
The Stirling number of the second kind ${n \brace k}$ (see \eqref{7}) enumerates the number of ways a set of $n$ elements can be partitioned into $k$ nonempty subsets, while the $n$th Bell number $\phi_{n}=\phi_{n}(1)$ (see \eqref{6}) counts the number of ways a set of $n$ elements can be partitioned into nonempty subsets.
The classical unsigned Lah number $L(n,k)$ (see \eqref{3}, \eqref{4}, \eqref{5}) counts the number of ways a set of $n$ elements can be partitioned into $k$ nonempty linearly ordered subsets. Just as the Bell number is to the Stirling numbers, the Lah-Bell number is to the unsigned Lah numbers. More precisely, the $n$th Lah-Bell number $B_{n}^{L}=B_{n}^{L}(1)$ (see \eqref{8}, \eqref{9}) is defined as the number of ways a set of $n$ elements can be partitioned into nonempty linearly ordered subsets. \par
The aim of this paper is to study probabilistic versions of the Lah numbers $L(n,k)$, and the Lah-Bell polynomials $B_{n}^{L}(x)$, which are a natural extension of Lah-Bell numbers, namely the probabilisitc Lah numbers associated with $Y$ and the probabilistic Lah- Bell polynomials associated with $Y$. Here $Y$ is a random variable whose moment generating function exists in some neighborhood of the origin. We derive some properties, explicit expressions, recurrence relations and certain identities for those numbers and polynomials. In addition, we treat the special cases that $Y$ is the Poisson random variable with parameter $\alpha >0$ and the Bernoulli random variable with probability of success $p$. \par
The outline of this paper is as follows. In Section 1, we recall the Stirling numbers of the second kind and the Bell polynomials. We remind the reader of the Lah numbers and the Lah-Bell polynomials. We recall the partial Bell polynomials and the complete Bell polynomials. Section 2 is the main result of this paper. Assume that $Y$ is a random variable such that the moment generating function of $Y$,\,\, $E[e^{tY}]=\sum_{n=0}^{\infty}\frac{t^{n}}{n!}E[Y^{n}], \quad (|t| <r)$, exists for some $r >0$. Let $(Y_{j})_{j\ge 1}$ be a sequence of mutually independent copies of the random variable $Y$, and let $S_{k}=Y_{1}+Y_{2}+\cdots+Y_{k},\,\, (k \ge 1)$,\,\, with $S_{0}=0$. Then we first define the probabilistic Lah numbers associated with the random variable $Y$, $L_{Y}(n,k)$, which are probabilistic versions of the Lah numbers and reduce to those numbers for $Y=1$. In Theorem 2.1, we derive an expression of $L_{Y}(n,k)$ in terms of $E\big[\langle S_{l}\rangle_{n}\big],\,\,(l=0,\dots,k)$. In Theorem 2.2, we show that $L_{Y}(n,k)= \frac{1}{k!}E\big[\triangle_{Y_{1},Y_{2},\dots,Y_{k}}\langle 0\rangle_{n}\big]$,\,\,(see \eqref{18}).
Next, we define the probabilistic Lah-Bell polynomials associated with the random variable $Y$, $B_{n}^{(L,Y)}(x)$, which are probabilistic versions of Lah-Bell polynomials. We obtain in Theorem 2.3 the generating function of the polynomials $B_{n}^{(L,Y)}(x)$, from which a Dobinski-like formula for those ones are deduced in Theorem 2.4. Two expressions of $B_{n}^{(L,Y)}(x)$ are obtained in terms of the partial Bell polynomilas in Theorems 2.5 and 2.8. A recurrence relation for $B_{n}^{(L,Y)}(x)$ is derived in Theorem 2.6. It is shown that $B_{n}^{(L,Y)}(x)$ satisfies the binomial identity in Theorem 2.7. Two identities involving $L_{Y}(n,k)$ and $B_{n}^{(L,Y)}(x)$ are obtained in Theorems 2.9 and 2.10. A finite sum expression for the $k$th derivative of $B_{n}^{(L,Y)}(x)$ is shown in Theorem 2.11. Let $Y$ be the Poisson random variable with parameter $\alpha >0$. Then we show $E\Big[\langle S_{k}\rangle_{n}\Big]=B_{n}^{L}(\alpha k)$ in Theorem 2.12 and in terms of the Lah numbers and the Bell polynomials a finite sum expression of $B_{n}^{(L,Y)}(x)$ in Theorem 2.13. Finally, in Theorem 2.14 we show $B_{n}^{(L,Y)}(x)=B_{n}^{L}(px)$ and $L_{Y}(n,k)=p^{k}L(n,k)$ when $Y$ is the Bernoulli random variable with probability of success $p$.\par

\vspace{0.1in}

We recall that the falling factorial sequence is defined by
\begin{equation}
(x)_{0}=1,\quad (x)_{n}=x(x-1)(x-2)\cdots(x-n+1),\quad (n\ge 1),\quad (\mathrm{see}\ [1-25]).\label{1}
\end{equation}
Also, the rising factorial sequence is given by
\begin{equation}
\langle x\rangle_{0}=1,\quad \langle x\rangle_{n}=x(x+1)(x+2)\cdots(x+n-1),\quad (n\ge 1).\label{2}
\end{equation}
For any integers $n,k $ with $n \ge k \ge 0$, the unsigned Lah numbers are defined by
\begin{equation}
\langle x\rangle_{n}=\sum_{k=0}^{n}L(n,k)(x)_{k},\quad (n\ge 0),\quad (\mathrm{see}\ [9,12,13,21]).\label{3}
\end{equation}
Thus, by \eqref{3}, we get
\begin{equation}
\frac{1}{k!}\bigg(\frac{1}{1-t}-1\bigg)^{k}=\sum_{n=k}^{\infty}L(n,k)\frac{t^{n}}{n!},\quad (k\ge 0),\quad (\mathrm{see}\ [11-13,15,16,18,19]).\label{4}
\end{equation}
From \eqref{4}, we note that
\begin{equation}
L(n,k)=\frac{n!}{k!}\binom{n-1}{k-1},\quad (n \ge k\ge 0),\quad (\mathrm{see}\ [9,12,13,22]).\label{5}
\end{equation}
It is well known that the Bell polynomials are defined by
\begin{equation}
\phi_{n}(x)=\sum_{k=0}^{n}{n \brace k}x^{k},\quad (\mathrm{see}\ [14,22]). \label{6}
\end{equation}
where ${n \brace k}$ are the Stirling numbers of the second kind given by
\begin{equation}
x^{n}=\sum_{k=0}^{n}{n\brace k}(x)_{k},\quad (\mathrm{see}\ [14,17,20-22]). \label{7}	
\end{equation}
In view of \eqref{6}, Lah-Bell polynomials are define by
\begin{equation}
B_{n}^{L}(x)=\sum_{l=0}^{n}L(n,l)x^{l},\quad (n\ge 0),\quad (\mathrm{see}\ [12,13]). \label{8}
\end{equation}
Thus, by \eqref{4} and \eqref{8}, we get
\begin{equation}
e^{x(\frac{1}{1-t}-1)}=\sum_{n=0}^{\infty}B_{n}^{L}(x)\frac{t^{n}}{n!},\quad (\mathrm{see}\ [11-20]).\label{9}
\end{equation}
For $k\ge 0$, the partial Bell polynomials are defined by
\begin{equation}
\frac{1}{k!}\bigg(\sum_{j=1}^{\infty}x_{j}\frac{t^{j}}{j!}\bigg)^{k}=\sum_{n=k}^{\infty}B_{n,k}(x_{1},x_{2},\dots,x_{n-k+1})\frac{t^{n}}{n!},\quad (\mathrm{see}\ [21]).\label{10}	
\end{equation}
Thus, by \eqref{10}, we get
\begin{equation}
\begin{aligned}
&B_{n,k}(x_{1},x_{2},\dots,x_{n-k+1})\\
&=\sum_{\substack{k_{1}+k_{2}+\cdots+k_{n-k+1}=k \\ k_{1}+2k_{2}+\cdots+(n-k+1)k_{n-k+1}=n}}	\frac{n!}{k_{1}! k_{2}!\cdots k_{n-k+1}!}\bigg(\frac{x_{1}}{1!}\bigg)^{k_{1}} \bigg(\frac{x_{2}}{2!}\bigg)^{k_{2}}\cdots \bigg(\frac{x_{n-k+1}}{(n-k+1)!}\bigg)^{k_{n-k+1}}.
\end{aligned}\label{11}
\end{equation}
The complete Bell polynomials are given by
\begin{equation}
\exp\bigg(\sum_{j=1}^{\infty}x_{j}\frac{t^{j}}{j!}\bigg)=\sum_{n=0}^{\infty}B_{n}(x_{1},x_{2},\dots,x_{n})\frac{t^{n}}{n!},\quad (\mathrm{see}\ [9,16,19]).\label{12}
\end{equation}
Thus, by \eqref{10} and \eqref{11}, we get
\begin{equation}
	B_{n}(x_{1},x_{2},\dots,x_{n})=\sum_{k=0}^{n}B_{n,k}(x_{1},x_{2},\dots,x_{n-k+1}),\quad (\mathrm{see}\ [9,13,16]). \label{13}
 \end{equation}

\section{Probabilistic Lah numbers and Lah-Bell polynomials}

{\it{Throughout this section, we assume that $Y$ is a random variable such that the moment generating function of $\displaystyle Y,\ E\big[e^{tY}\big]=\sum_{n=0}^{\infty}E[Y^{n}]\frac{t^{n}}{n!}\displaystyle$, $(|t|<r)$, exists for some $r>0$}},\,\,(see [3,16]). {\it{We let $(Y_{j})_{j\ge 1}$ be a sequence of mutually independent copies of the random variable $Y$, and let}}
\begin{equation}
S_{0}=0,\quad S_{k}=Y_{1}+Y_{2}+\cdots+Y_{k},\quad (k\in\mathbb{N}),\,\,(\mathrm{see}\,\, [3,16]).\label{15}	
\end{equation}

We define the {\it{probabilistic Lah numbers asssociated with the random variable $Y$}} by
\begin{equation}
\frac{1}{k!}\Bigg(E\bigg[\bigg(\frac{1}{1-t}\bigg)^{Y}\bigg]-1\Bigg)^{k}=\sum_{n=k}^{\infty}L_{Y}(n,k)\frac{t^{n}}{n!},\quad (k\ge 0).\label{16}
\end{equation}
When $Y=1$, we note that $L_{Y}(n,k)=L(n,k)$. From \eqref{16}, we note that
\begin{align}
&\sum_{n=k}^{\infty}L_{Y}(n,k)\frac{t^{n}}{n!}=\frac{1}{k!}\sum_{l=0}^{k}\binom{k}{l}(-1)^{k-l}E\Bigg[\bigg(\frac{1}{1-t}\bigg)^{Y_{1}+Y_{2}+\cdots+Y_{l}}\Bigg] \label{17}\\
&=\frac{1}{k!}\sum_{l=0}^{k}\binom{k}{l}(-1)^{k-l}E\Bigg[\bigg(\frac{1}{1-t}\bigg)^{S_{l}}\Bigg]=\sum_{n=0}^{\infty}\frac{1}{k!}\sum_{l=0}^{k}\binom{k}{l}(-1)^{k-l}E\Big[\langle S_{l}\rangle_{n}\Big]\frac{t^{n}}{n!}.\nonumber	
\end{align}
Therefore, by comparing the coefficients on both sides of \eqref{17}, we obtain the following theorem.
\begin{theorem}
For any integers $n,k$ with $n\ge k\ge 0$, we have
\begin{displaymath}
L_{Y}(n,k)= \frac{1}{k!}\sum_{l=0}^{k}\binom{k}{l}(-1)^{k-l}E\big[\langle S_{l}\rangle_{n}\big].
\end{displaymath}
\end{theorem}
The operators $\triangle_{y}$ and $\triangle_{y_{1},y_{2},\dots,y_{m}}$ are respectively given by
\begin{equation}
\triangle_{y}f(x)=f(x+y)-f(x),\quad \triangle_{y_{1},y_{2},\dots,y_{m}}f(x)= \triangle_{y_{1}}\circ \triangle_{y_{2}}\circ \triangle_{y_{3}}\circ\cdots\circ \triangle_{y_{m}}f(x), \label{18}
\end{equation}
where $y \in \mathbb{R},\,\, and \,\, (y_{1},y_{2},\dots,y_{m})\in\mathbb{R}^{m}.$
Then, by \eqref{18}, we get
\begin{equation}
\begin{aligned}
\triangle_{y_{1},y_{2},\dots,y_{m}}\bigg(\frac{1}{1-t}\bigg)^{x}&= \bigg(\bigg(\frac{1}{1-t}\bigg)^{y_{1}}-1\bigg) \bigg(\bigg(\frac{1}{1-t}\bigg)^{y_{2}}-1\bigg)\cdots \bigg(\bigg(\frac{1}{1-t}\bigg)^{y_{m}}-1\bigg)\bigg(\frac{1}{1-t}\bigg)^{x}\\
&=\sum_{n=m}^{\infty}\frac{t^{n}}{n!} \triangle_{y_{1},y_{2},\dots,y_{m}}\langle x\rangle_{n}.
\end{aligned}\label{19}
\end{equation}
From \eqref{16}, \eqref{19}, we note that
\begin{align}
&\frac{1}{m!}\bigg(E\bigg[\bigg(\frac{1}{1-t}\bigg)^{Y}\bigg]-1\bigg)^{m}=\frac{1}{m!}\bigg(E\bigg[\bigg(\frac{1}{1-t}\bigg)^{Y}-1\bigg]\bigg)^{m}\label{20}\\
&=\frac{1}{m!}\bigg(E\bigg[\bigg(\bigg(\frac{1}{1-t}\bigg)^{Y_{1}}-1\bigg) \bigg(\bigg(\frac{1}{1-t}\bigg)^{Y_{2}}-1\bigg)\cdots \bigg(\bigg(\frac{1}{1-t}\bigg)^{Y_{m}}-1\bigg)\bigg]\nonumber\\
&=\sum_{n=m}^{\infty}\frac{1}{m!}E\big[\triangle_{Y_{1},Y_{2},\dots,Y_{m}}\langle 0\rangle_{n}\big]\frac{t^{n}}{n!}. \nonumber
\end{align}
Therefore, by \eqref{20}, we obtain the following theorem.
\begin{theorem}
For any integers $n,m$ with $n\ge m\ge 0$, we have
\begin{displaymath}
L_{Y}(n,m)= \frac{1}{m!}E\big[\triangle_{Y_{1},Y_{2},\dots,Y_{m}}\langle 0\rangle_{n}\big].
\end{displaymath}
\end{theorem}
We define the {\it{probabilistic Lah-Bell polynomials associated with the random variable $Y$}} by
\begin{equation}
B_{n}^{(L,Y)}(x)=\sum_{k=0}^{n}L_{Y}(n,k)x^{k},\quad (n\ge 0).\label{21}	
\end{equation}
When $Y=1$, we have $B_{n}^{(L,Y)}(x)=B_{n}^{L}(x)$,\quad (see \eqref{8}). \par
Form \eqref{21}, we derive the generating function of the probabilistic Lah-Bell polynomials associated the random variable $Y$, which is given by
\begin{align}
&\sum_{n=0}^{\infty}B_{n}^{(L,Y)}(x)\frac{t^{n}}{n!}=\sum_{n=0}^{\infty}\bigg(\sum_{k=0}^{n}L_{Y}(n,k)x^{k}\bigg)\frac{t^{n}}{n!}. \label{22}\\
&=\sum_{k=0}^{\infty}x^{k}\bigg(\sum_{n=k}^{\infty}L_{Y}(n,k)\frac{t^{n}}{n!}\bigg)=\sum_{k=0}^{\infty}x^{k}\bigg(E\bigg[\bigg(\frac{1}{1-t}\bigg)^{Y}\bigg]-1\bigg)^{k}\frac{1}{k!}\nonumber\\
&=\exp\bigg(x\bigg(E\bigg[\bigg(\frac{1}{1-t}\bigg)^{Y}\bigg]-1\bigg)\bigg).\nonumber
\end{align}
Thus, by \eqref{22}, we obtain the following theorem.
\begin{theorem}
The generating function of the probabilistic Lah-Bell polynomials of the random variable $Y$ is given by
\begin{equation}
e^{x\big(E\big[\big(\frac{1}{1-t}\big)^{Y}\big]-1\big)}=\sum_{n=0}^{\infty}B_{n}^{(L,Y)}(x)\frac{t^{n}}{n!}	.\label{23}
\end{equation}
\end{theorem}
\noindent By \eqref{23}, we get
\begin{align}
\sum_{n=0}^{\infty}B_{n}^{(L,Y)}(x)\frac{t^{n}}{n!}&= e^{x\big(E\big[\big(\frac{1}{1-t}\big)^{Y}\big]-1\big)}=e^{-x}e^{xE\big[\big(\frac{1}{1-t}\big)^{Y}\big]}\label{24} \\
&=e^{-x}\sum_{k=0}^{\infty}\frac{x^{k}}{k!}\bigg(E\bigg[\bigg(\frac{1}{1-t}\bigg)^{Y}\bigg]\bigg)^{k}=e^{-x}\sum_{k=0}^{\infty}\frac{x^{k}}{k!} E\bigg[\bigg(\frac{1}{1-t}\bigg)^{Y_{1}}\cdots \bigg(\frac{1}{1-t}\bigg)^{Y_{k}}\bigg]\nonumber\\
&=e^{-x}\sum_{k=0}^{\infty}\frac{x^{k}}{k!}E\bigg[\bigg(\frac{1}{1-t}\bigg)^{S_{k}}\bigg]=\sum_{n=0}^{\infty}\bigg(e^{-x}\sum_{k=0}^{\infty}\frac{x^{k}}{k!}E\big[\langle S_{k}\rangle_{n}\big]\bigg)\frac{t^{n}}{n!}.\nonumber
\end{align}
Therefore, by \eqref{24}, we obtain the following Dobinski-like formula.
\begin{theorem}
For any integer $n\ge 0$, we have
\begin{displaymath}
B_{n}^{(L,Y)}(x)=e^{-x}\sum_{k=0}^{\infty}\frac{x^{k}}{k!}E\big[\langle S_{k}\rangle_{n}\big].
\end{displaymath}
\end{theorem}
From \eqref{23}, we note that
\begin{align}
&\sum_{n=0}^{\infty}B_{n}^{(L,Y)}(x)\frac{t^{n}}{n!}=e^{x\big(E\big[\big(\frac{1}{1-t}\big)^{Y}\big]-1\big)}\label{25} \\
&=e^{x\sum_{j=1}^{\infty}E[\langle Y\rangle_{j}]\frac{t^{j}}{j!}}=\sum_{k=0}^{\infty}\frac{1}{k!}\bigg(x\sum_{j=1}^{\infty}E[\langle Y\rangle_{j}]\frac{t^{j}}{j!}\bigg)^{k}\nonumber\\
&=\sum_{k=0}^{\infty}\sum_{n=k}^{\infty}B_{n,k}\big(xE[\langle Y\rangle_{1}], xE[\langle Y\rangle_{2}],\dots, xE[\langle Y\rangle_{n-k+1}]\big)\frac{t^{n}}{n!}\nonumber\\
&=\sum_{n=0}^{\infty}\bigg(\sum_{k=0}^{n}B_{n,k}\big(xE[\langle Y\rangle_{1}], xE[\langle Y\rangle_{2}],\dots, xE[\langle Y\rangle_{n-k+1}]\big)\bigg)\frac{t^{n}}{n!}.\nonumber
\end{align}
Therefore, by comparing the coefficients on both sides of \eqref{25}, we obtain the following theorem.
\begin{theorem}
For any integers $n,k$ with $n\ge k\ge 0$, we have
\begin{displaymath}
B_{n}^{(L,Y)}(x)=\sum_{k=0}^{n} B_{n,k}\big(xE[\langle Y\rangle_{1}], xE[\langle Y\rangle_{2}],\dots, xE[\langle Y\rangle_{n-k+1}]\big).
\end{displaymath}
\end{theorem}
Now, we observe that
\begin{align}
&\sum_{n=0}^{\infty}B_{n+1}^{(L,Y)}(x)\frac{t^{n}}{n!}=\frac{d}{dt}\sum_{n=0}^{\infty}B_{n}^{(L,Y)}(x)\frac{t^{n}}{n!}=\frac{d}{dt}e^{x\big(E\big[\big(\frac{1}{1-t}\big)^{Y}\big]-1\big)}\label{26} \\
&=xE\bigg[Y\bigg(\frac{1}{1-t}\bigg)^{Y+1}\bigg] e^{x\big(E\big[\big(\frac{1}{1-t}\big)^{Y}\big]-1\big)}=x\sum_{k=0}^{\infty}E\Big[\langle Y\rangle_{k+1}\Big]\frac{t^{k}}{k!}\sum_{m=0}^{\infty}B_{m}^{(L,Y)}(x)\frac{t^{m}}{m!} \nonumber \\
&=\sum_{n=0}^{\infty}\bigg(x\sum_{k=0}^{n}\binom{n}{k}E\Big[\langle Y\rangle_{k+1}\Big]B_{n-k}^{(L,Y)}(x)\bigg)\frac{t^{n}}{n!}. \nonumber	
\end{align}
Therefore, by \eqref{26}, we obtain the following theorem.
\begin{theorem}
For any integer $n\ge 0$, we have the recurrence relation
\begin{displaymath}
B_{n+1}^{(L,Y)}(x)=x\sum_{k=0}^{n}\binom{n}{k}E\Big[\langle Y\rangle_{k+1}\Big]B_{n-k}^{(L,Y)}(x).
\end{displaymath}
\end{theorem}
From \eqref{23}, we have
\begin{align}
&\sum_{n=0}^{\infty}B_{n}^{(L,Y)}(x+y)\frac{t^{n}}{n!}=e^{(x+y)\big(E\big[\big(\frac{1}{1-t}\big)^{Y}\big]-1\big)}=e^{x\big(E\big[\big(\frac{1}{1-t}\big)^{Y}-1\big]\big)}e^{y\big(E\big[\big(\frac{1}{1-t}\big)^{Y}-1\big]\big)}\label{27}\\
&=\sum_{k=0}^{\infty}B_{k}^{(L,Y)}(x)\frac{t^{k}}{k!}\sum_{m=0}^{\infty}B_{m}^{(L,Y)}(y)\frac{t^{m}}{m!}=\sum_{n=0}^{\infty}\bigg(\sum_{k=0}^{n}\binom{n}{k}B_{k}^{(L,Y)}(x)B_{n-k}^{(L,Y)}(y)\bigg)\frac{t^{n}}{n!}.\nonumber
\end{align}
Therefore, by \eqref{27}, we obtain the following binomial identity.
\begin{theorem}
For any integer $n\ge 0$, we have
\begin{displaymath}
B_{n}^{(L,Y)}(x+y)=\sum_{k=0}^{n}\binom{n}{k}B_{k}^{(L,Y)}(x)B_{n-k}^{(L,Y)}(y).
\end{displaymath}
\end{theorem}
By \eqref{23}, we see that
\begin{align}
&\sum_{n=0}^{\infty}B_{n}^{(L,Y)}(x)\frac{t^{n}}{n!}= e^{x\big(E\big[\big(\frac{1}{1-t}\big)^{Y}\big]-1\big)} \label{28} \\
&=\sum_{k=0}^{\infty}\binom{x}{k}\Big(e^{\big(E\big[\big(\frac{1}{1-t}\big)^{Y}\big]-1\big)}-1\Big)^{k}=\sum_{k=0}^{\infty}k!\binom{x}{k}\frac{1}{k!}\bigg(\sum_{j=1}^{\infty}B_{j}^{(L,Y)}(1)\frac{t^{j}}{j!}\bigg)^{k}\nonumber\\
&=\sum_{k=0}^{\infty}k!\binom{x}{k}\sum_{n=k}^{\infty}B_{n,k}\Big(B_{1}^{(L,Y)}(1),B_{2}^{(L,Y)}(1),\dots,B_{n-k+1}^{(L,Y)}(1)\Big)\frac{t^{n}}{n!}\nonumber \\
&=\sum_{n=0}^{\infty}\bigg(\sum_{k=0}^{n}k!\binom{x}{k} B_{n,k}\Big(B_{1}^{(L,Y)}(1),B_{2}^{(L,Y)}(1),\dots,B_{n-k+1}^{(L,Y)}(1)\Big)\Bigg)\frac{t^{n}}{n!}.\nonumber
\end{align}
Therefore, by \eqref{28}, we obtain the following theroem.
\begin{theorem}
For any integer $n\ge 0$, we ahve
\begin{displaymath}
B_{n}^{(L,Y)}(x)=\sum_{k=0}^{n}k!\binom{x}{k} B_{n,k}\Big(B_{1}^{(L,Y)}(1),B_{2}^{(L,Y)}(1),\dots,B_{n-k+1}^{(L,Y)}(1)\Big).
\end{displaymath}
\end{theorem}
Note that
\begin{equation}
te^{x\big(E\big[\big(\frac{1}{1-t}\big)^{Y}\big]-1\big)}=t\sum_{j=0}^{\infty}B_{j}^{(L,Y)}(x)\frac{t^{j}}{j!}=\sum_{j=1}^{\infty}jB_{j-1}^{(L,Y)}(x)\frac{t^{j}}{j!}.\label{29}	
\end{equation}
By \eqref{29} and for $k \ge 0$, we get
\begin{align}
&\bigg(\sum_{j=1}^{\infty}jB_{j-1}^{(L,Y)}(x)\frac{t^{j}}{j!}\bigg)^{k}=t^{k}e^{kx\big(E\big[\big(\frac{1}{1-t}\big)^{Y}\big]-1\big)}\label{30}\\
&=t^{k}\sum_{j=0}^{\infty}k^{j}x^{j}\frac{1}{j!}\bigg(E\bigg[\bigg(\frac{1}{1-t}\bigg)^{Y} \bigg]-1\bigg)^{j}=t^{k}\sum_{j=0}^{\infty}k^{j}x^{j}\sum_{n=j}^{\infty}L_{Y}(n,j)\frac{t^{n}}{n!}\nonumber\\
&=\sum_{n=0}^{\infty}\sum_{j=0}^{n}k^{j}x^{j}L_{Y}(n,j)\frac{t^{n+k}}{n!}=\sum_{n=k}^{\infty}\sum_{j=0}^{n-k}k^{j}x^{j}L_{Y}(n-k,j)\frac{t^{n}}{(n-k)!} \nonumber \\
&=\sum_{n=k}^{\infty}k!\sum_{j=0}^{n-k}\binom{n}{k}k^{j}x^{j}L_{Y}(n-k,j)\frac{t^{n}}{n!}. \nonumber
\end{align}
From \eqref{30}, we have
\begin{align}
&\sum_{n=k}^{\infty}\sum_{j=0}^{n-k}\binom{n}{k}k^{j}x^{j}L_{Y}(n-k,j)\frac{t^{n}}{n!}=\frac{1}{k!}\bigg(\sum_{j=1}^{\infty}jB_{j-1}^{(L,Y)}(x)\frac{t^{j}}{j!}\bigg)^{k} \label{31} \\
&=\sum_{n=k}^{\infty}B_{n,k}\Big(B_{0}^{(L,Y)}(x),2B_{1}^{(L,Y)}(x),\dots,(n-k+1)B_{n-k}^{(L,Y)}(x)\Big)\frac{t^{n}}{n!}.\nonumber
\end{align}
Therefore, by \eqref{31}, we obtain the following theorem.
\begin{theorem}
For any integers $n,k$ with $n\ge k\ge 0$, we have
\begin{align*}
&\sum_{j=0}^{n-k}\binom{n}{k}k^{j}x^{j}L_{Y}(n-k,j)\\
&= B_{n,k}\Big(B_{0}^{(L,Y)}(x),2B_{1}^{(L,Y)}(x),\dots,(n-k+1)B_{n-k}^{(L,Y)}(x)\Big).
\end{align*}
\end{theorem}
Now, we observe that
\begin{align}
&\sum_{n=k}^{\infty} B_{n,k}\Big(B_{1}^{(L,Y)}(x),B_{2}^{(L,Y)}(x),\dots,B_{n-k+1}^{(L,Y)}(x)\Big)\frac{t^{n}}{n!}\label{32} \\
&=\frac{1}{k!}\bigg(\sum_{j=1}^{\infty}B_{j}^{(L,Y)}(x)\frac{t^{j}}{j!}\bigg)^{k}=\frac{1}{k!}\bigg(e^{x\big(E\big[\big(\frac{1}{1-t}\big)^{Y}\big]-1\big)}-1\bigg)^{k} \nonumber\\
&=\sum_{j=k}^{\infty}{j \brace k}x^{j}\frac{1}{j!}\bigg(E\bigg[\bigg(\frac{1}{1-t}\bigg)^{Y}\bigg]-1\bigg)^{j}\nonumber \\
&=\sum_{j=k}^{\infty}{j \brace k}x^{j}\sum_{n=j}^{\infty}L_{Y}(n,j)\frac{t^{n}}{n!}=\sum_{n=k}^{\infty}\bigg(\sum_{j=k}^{n}{j \brace k}L_{Y}(n,j)x^{j}\bigg)\frac{t^{n}}{n!}. \nonumber
\end{align}
Therefore, by comparing the coefficients on both sides of \eqref{32}, we obtain following theorem.
\begin{theorem}
For any integers $n,k$ with $n\ge k\ge 0$, we have
\begin{displaymath}
B_{n,k}\Big(B_{1}^{(L,Y)}(x),B_{2}^{(L,Y)}(x),\dots,B_{n-k+1}^{(L,Y)}(x)\Big)=\sum_{j=k}^{n}{j\brace k}L_{Y}(n,j)x^{j}.
\end{displaymath}
\end{theorem}
From \eqref{23} and $k \ge 0$, we have
\begin{align}
&\sum_{n=0}^{\infty}\bigg(\frac{d}{dx}\bigg)^{k}B_{n}^{(L,Y)}(x)\frac{t^{n}}{n!}=\bigg(\frac{d}{dx}\bigg)^{k}e^{x\big(E\big[\big(\frac{1}{1-t}\big)^{Y}\big]-1\big)}\label{33} \\
&=k!\frac{1}{k!}\bigg(E\bigg[\bigg(\frac{1}{1-t}\bigg)^{Y}\bigg]-1\bigg)^{k} e^{x\big(E\big[\big(\frac{1}{1-t}\big)^{Y}\big]-1\big)}\nonumber\\
&=k!\sum_{l=k}^{\infty}L_{Y}(l,k)\frac{t^{l}}{l!}\sum_{j=0}^{\infty}B_{j}^{(L,Y)}(x)\frac{t^{j}}{j!}=\sum_{n=k}^{\infty}\bigg(k!\sum_{j=0}^{n-k}B_{j}^{(L,Y)}(x)L_{Y}(n-j,k)\binom{n}{j}\bigg)\frac{t^{n}}{n!}.\nonumber
\end{align}
In particular, for $k=1$, we get
\begin{align}
\sum_{n=1}^{\infty}\frac{d}{dx}B_{n}^{(L,Y)}(x)\frac{t^{n}}{n!}&=\bigg(E\bigg[\bigg(\frac{1}{1-t}\bigg)^{Y}\bigg]-1\bigg) e^{x\big(E\big[\big(\frac{1}{1-t}\big)^{Y}\big]-1\big)} \label{34} \\
&=\sum_{l=1}^{\infty}E\big[\langle Y\rangle_{l}\big]\frac{t^{l}}{l!}\sum_{j=0}^{\infty}B_{j}^{(L,Y)}(x)\frac{t^{j}}{j!}\nonumber \\
&=\sum_{n=1}^{\infty}\bigg(\sum_{j=0}^{n-1}E\big[\langle Y\rangle_{n-j}\big]B_{j}^{(L,Y)}(x)\binom{n}{j}\bigg)\frac{t^{n}}{n!}. \nonumber 	
\end{align}
Therefore, by \eqref{33} and \eqref{34}, we obtain the following theorem.
\begin{theorem}
For any integers $n,k$ with $n\ge k \ge 0$, we have
\begin{displaymath}
\bigg(\frac{d}{dx}\bigg)^{k}B_{n}^{(L,Y)}(x)=k!\sum_{j=0}^{n-k}\binom{n}{j}B_{j}^{(L,Y)}(x)L_{Y}(n-j,k).
\end{displaymath}
In particular, for $k=1$, we get
\begin{displaymath}
\frac{d}{dx}B_{n}^{(L,Y)}(x)=\sum_{j=0}^{n-1}\binom{n}{j}E\big[\langle Y\rangle_{n-j}\big]B_{j}^{(L,Y)}(x),\quad (n \ge 1).
\end{displaymath}
\end{theorem}
\noindent From \eqref{16}, we have
\begin{align}
&\sum_{n=k}^{\infty}L_{Y}(n,k)\frac{t^{n}}{n!}=\frac{1}{k!}\bigg(E\bigg[\bigg(\frac{1}{1-t}\bigg)^{Y}\bigg]-1\bigg)^{k}=\frac{1}{k!}\bigg(\sum_{j=1}^{\infty}E\big[\langle Y\rangle_{j}\big]\frac{t^{j}}{j!}\bigg)^{k}\label{35} \\
&=\sum_{n=k}^{\infty}B_{n,k}\bigg(E\Big[\langle Y\rangle_{1}\Big], E\big[\langle Y\rangle_{2}\Big],\dots, E\big[\langle Y\rangle_{n-k+1}\Big]\bigg)\frac{t^{n}}{n!}.\nonumber
\end{align}
Let $Y$ be the Poission random variable with parameter $\alpha>0$. Then we have
\begin{align}
&\sum_{n=0}^{\infty}E\Big[\langle S_{k}\rangle_{n}\Big]\frac{t^{n}}{n!}=E\bigg[\bigg(\frac{1}{1-t}\bigg)^{S_{k}}\bigg]=E\bigg[\bigg(\frac{1}{1-t}\bigg)^{Y_{1}+\cdots+Y_{k}}\bigg] \label{36} \\
&=E\bigg[\bigg(\frac{1}{1-t}\bigg)^{Y_{1}}\bigg] E\bigg[\bigg(\frac{1}{1-t}\bigg)^{Y_{2}}\bigg]\cdots E\bigg[\bigg(\frac{1}{1-t}\bigg)^{Y_{k}}\bigg]\nonumber \\
&=e^{-\alpha k}e^{\alpha k\big(\frac{1}{1-t}\big)}=e^{\alpha k\big(\frac{1}{1-t}-1\big)}=\sum_{n=0}^{\infty}B_{n}^{L}(\alpha k)\frac{t^{n}}{n!}.\nonumber
\end{align}
Therefore, by \eqref{36}, we obtain the following theorem.
\begin{theorem}
Let $Y$ be the Poisson random variable with parameter $\alpha>0$. Then we have
\begin{displaymath}
E\Big[\langle S_{k}\rangle_{n}\Big]=B_{n}^{L}(\alpha k),\quad (n, k\ge 0).
\end{displaymath}
\end{theorem}
From \eqref{23} and \eqref{27}, we note that
\begin{align}
\sum_{n=0}^{\infty}B_{n}^{(L,Y)}(x)\frac{t^{n}}{n!}&=e^{x\big(E\big[\big(\frac{1}{1-t}\big)^{Y}\big]-1\big)}=e^{x(e^{\alpha(\frac{1}{1-t}-1)}-1)}\label{38}\\
&=\sum_{k=0}^{\infty}\phi_{k}(x)\alpha^{k}\frac{1}{k!}\bigg(\frac{1}{1-t}-1\bigg)^{k}=\sum_{k=0}^{\infty}\phi_{k}(x)\alpha^{k}\sum_{n=k}^{\infty}L(n,k)\frac{t^{n}}{n!}\nonumber\\
&=\sum_{n=0}^{\infty}\bigg(\sum_{k=0}^{n}\phi_{k}(x)\alpha^{k}L(n,k)\bigg)\frac{t^{n}}{n!}. \nonumber
\end{align}
Therefore, by \eqref{38}, we obtain the following theorem.
\begin{theorem}
Let $Y$ be the Poisson random variable with parameter $\alpha>0$. For $n\ge 0$, we have
\begin{displaymath}
B_{n}^{(L,Y)}(x)=\sum_{k=0}^{n}\phi_{k}(x)\alpha^{k}L(n,k).
\end{displaymath}
\end{theorem}
Assume that $Y$ is the Bernoulli random variable with probability of success $p$. Then, we have
\begin{equation}
E\big[\langle Y\rangle_{n}\big]=\sum_{i=0}^{1}\langle i\rangle_{n}p(i)=n!p,\quad(n \ge 1).\label{39}	
\end{equation}
Thus, by \eqref{35} and \eqref{39}, we get
\begin{align}
L_{Y}(n,k)&=B_{n,k}\big(1!p,2!p,\dots,(n-k+1)!p\big)\label{40}\\	
&=p^{k}B_{n,k}\big(1!,2!,\dots,(n-k+1)!\big)=p^{k}L(n,k),\quad(n \ge k,\, n \ge 1),\nonumber
\end{align}
and, by \eqref{8} and \eqref{21}, we have
\begin{equation}
B_{n}^{(L,Y)}(x)=\sum_{k=0}^{n}L_{Y}(n,k)x^{k}=\sum_{k=0}^{n}(px)^{k}L(n,k)=B_{n}^{L}(px),\quad (n\ge 1). \label{41}
\end{equation}
Therefore, by \eqref{41}, we obtain the following theorem.
\begin{theorem}
Let $Y$ be the Bernoulli random variable with probability of success $p$. Then we have
\begin{displaymath}
B_{n}^{(L,Y)}(x)=B_{n}^{L}(px),\quad (n\ge 1).
\end{displaymath}
In addition, we have
\begin{displaymath}
L_{Y}(n,k)=p^{k}L(n,k),\quad (n \ge k,\, n \ge 1).
\end{displaymath}
\end{theorem}
\section{Conclusion}
In this paper, by using generating functions we studied the probabilistic Lah numbers associated with $Y$ and the probabilistic Lah-Bell polynomials associated with $Y$, as probabilistic versions of the Lah numbers and the Lah-Bell polynomials, respectively. Here $Y$ is a random variable such that the moment generating function of $Y$ exists in a neighborhood of the origin. In more detail, we derived some expressions for $L_{Y}(n,k)$ (see Theorems 2.1, 2.2), several expressions for $B_{n}^{(L,Y)}(x)$ (see Theorems 2.4, 2.5, 2.8), the generating function of $B_{n}^{(L,Y)}(x)$ (see Theorem 2.3), a recurrence relation for $B_{n}^{(L,Y)}(x)$ (see Theorem 2.6), the binomial identity for $B_{n}^{(L,Y)}(x)$ (see Theorem 2.7), and certain identities involving $L_{Y}(n,k)$, $B_{k}^{(L,Y)}(x)$ and the partial Bell polynomial (see Theorems 2.9, 2.10), and the higher order derivatives of $B_{n}^{(L,Y)}(x)$ (see Theorem 2.11). In addition, we obtained the identity $E\Big[\langle S_{k}\rangle_{n}\Big]=B_{n}^{L}(\alpha k)$ (see Theorem 2.12) and an expression of $B_{n}^{(L,Y)}(x)$ in terms of Bell polynomials and Lah numbers (see Theorem 2.13) when $Y$ is the Poisson random variable with parameter $\alpha >0$. Finally, we showed $B_{n}^{(L,Y)}(x)=B_{n}^{L}(px)$ and $L_{Y}(n,k)=p^{k}L(n,k)$ when $Y$ is the Bernoulli random variable with probability of success $p$ (see Theorem 2.14).\par
As one of our future projects, we would like to continue to investigate degenerate versions, $\lambda$-analogues and probabilistic versions of many special polynomials and numbers and to find their applications to physics, science and engineering as well as to mathematics.

\end{document}